\documentclass[amsart.cls,11pt]{article}
\newtheorem{lemma}{Lemma}
\newtheorem{definition}{Definition}
\newtheorem{theorem}{Theorem}
\newtheorem{conjecture}{Conjecture}
\newtheorem{remark}{Remark}
\newtheorem{corollary}{Corollary}
\begin{document}

\title{{\bf On the nonasymptotic prime number
distribution}\footnote{{\bf 1991 Mathematics Subject
Classification:} Primary 11N05, 11N35, 11N3b; Secondary 04A10.  {\it Key
Words and Phrases:} multiple Eratosthenes sieve, Eratosthenes progressions,
prime number separating theorem, matrices ${}^2 N$ and ${}^2 P$,
prime number inner distribution law (PNIDL), analytic form of PNIDL,
coagulates of prime numbers, strongly related twins.}}
\author{Lubomir Alexandrov  \\ Institute for Nuclear Research and Nuclear
Energy \\ 72 Tzarigrad Road, 1784 Sofia, Bulgaria} \maketitle

\begin{abstract}
The objective of this paper is to introduce an approach to the study of the
nonasymptotic distribution of prime numbers. The natural numbers are
represented by theorem 1 in the matrix form ${}^2 N$. The first column of
the infinite matrix ${}^2 N$ starts with the unit and contains all composite
numbers in ascending order. The infinite rows of this matrix except for the
first elements contain prime numbers only, which are determined by an
uniform recurrence law. At least one of the elements of the twin pairs of
prime numbers is an element of the second column of the matrix ${}^2 N$
(theorem 3). The basic information on the nonasymptotic prime number
distribution is contained in the distribution of the elements of the second
column of the matrix ${}^2 N$.
\end{abstract}

\section{Introduction}

The multiplicative and the additive structure of natural numbers is based on
prime numbers. The derivation of results on the distribution of prime
numbers in the set of natural numbers can affect almost all mathematical
theories and their applications. The objective of this paper is to introduce
an approach to the study of the nonasymptotic distribution of prime numbers.
The results obtained could be applicable in quantum physics, quantum
chemistry and molecular biology.

The Euclid theorem on the existence of an infinite number of prime numbers
\begin{equation}
\pi (x)\rightarrow \infty \ \ \ \ \ as\ x\rightarrow \infty ,
\end{equation}
(where $\pi (x)$ is the number of primes not exceeding $x$), the
Eratosthenes sieve formula
\begin{equation}
\pi (x)=\pi (\sqrt{x})-1+\sum_d(-1)^{\nu (d)}[\frac xd],
\end{equation}
(where $d$ runs over the divisors of the products of all primes not
exceeding $\sqrt{x}$, $\nu (d)$ is the number of the prime divisors of $d$
and $[u]$ is the integer part of $u$) as well as the asymptotic distribution
law of prime numbers (see f.i. [1] - [4])
\begin{equation}
\pi (x)=li\ x+O(xe^{-c\sqrt{\ln x}}) \ \ \ as\ x\rightarrow \infty ,\
c=const>1,
\end{equation}
where
\begin{equation}
li\
x=\int_2^x\frac{dt}{\ln t}=\frac x{\ln x}+\frac{1!x}{\ln ^2x}+\ldots +
\frac{(k-1)!x}{\ln ^kx}+O(\frac x{\ln ^{k-1}x})
\end{equation}
do not answer the question how often the prime numbers are encountered and
how they are distributed amidst the natural numbers when $x<\infty $. Even
if the Riemann fifth hypothesis [5] stating that the nontrivial solutions of
the equation
\begin{equation}
\zeta (s)\equiv \sum_{n=1}^\infty \frac 1{n^s}=0,\ \ s=\sigma +it
\end{equation}
lie on the complex straight line $(\sigma =1/2,\ t)$ was proved it could not
give an answer to that question too. The reason is that the fourth Riemann
hypothesis [5] (proved by H. von Mangoldt) stating that the expression
\[
\ P_0(x)=li\ x-\sum_\rho li\ x^\rho +\int_x^\infty \frac{du}{(u^2-1)\ln u}
-\ln 2
\]
where $\rho $ runs over all nontrivial solutions of eq. (5) and
\[
\ P_0(x)={\frac 12}(P(x+0)+P(x-0)),\ \ \ P(x)=\sum_{n\leq x}{\frac{\Lambda
(n)}{\ln x}}\ \ ,
\]
\[
\ \Lambda (n)=\left\{
\begin{array}{l}
\ln \rho ,\ \ when\ n=p^m\ \ \ and\ p\ is\ a\ prime,\ \ m\geq 1 \\
0\ \ \ otherwise
\end{array}
\right.
\]
only gives a connection between the distribution of the nontrivial zeroes of
eq. (5) and the problem of the prime numbers distribution.

In the present paper an approach based on the separation of subsets of prime
numbers for which explicit distribution law exists is used in the search for
an answer to the question is there any exact law for $\pi(x), \ 2 \leq x <
\infty$ or the prime numbers are spread amidst the natural numbers in total
disorder. The separation of these subsets is made by means of a new sieve
("the multiple Eratosthenes sieve" [6]). The sieve and its generalizations
will be published in a separate paper. In the presentation of the results on
the nonasymptotic distribution of prime numbers the sieve is present only
implicitly which simplifies the text. In the present paper all results are
obtained by elementary methods except the use of the inequality from [7].

\section{The prime number inner distribution law}

Let the set of prime numbers is enumerated in incremental order of elements
according to (1)

\[
P = \{2,3,5,7,11,13, \ldots \} = \{p_n\}_{n \in N} ,
\]
where $N$ is the set of natural numbers. Thus two reciprocal
number-theoretic functions are introduced
\[
\pi(p) : P \rightarrow N \ \ \ (i.e. \ \ \pi(p_n) = n)
\]
"the number $\pi(p)$ of the prime number $p$" and
\[
\pi^{-1}(n) : N \rightarrow P \ \ \ (i.e. \ \ \pi^{-1}(n) = p_n)
\]
"the prime number $\pi^{-1}(n)$ of number $n$".

It is evident that $\pi (p)$ and $\pi ^{-1}(n)$ satisfy the identities
\begin{equation}
\pi (\pi ^{-1}(n))=n
\end{equation}
and
\[
\ \pi ^{-1}(\pi (p))=p .
\]
These functions are also strictly monotonous
\[
\ p^{\prime }<p^{\prime \prime }\Longrightarrow \pi (p^{\prime })<\pi
(p^{\prime \prime }),\ \ p^{\prime },p^{\prime \prime }\in P,
\]
\begin{equation}
n^{\prime }<n^{\prime \prime }\Longrightarrow \pi ^{-1}(n^{\prime })<\pi
^{-1}(n^{\prime \prime }),\ \ n^{\prime },n^{\prime \prime }\in N.
\end{equation}
The following auxiliary proposition is a corollary from the inequality $
p_n>n\ln n$ from [7].

\begin{lemma}
For any prime number $n\geq 2$ the following estimates are correct:
\begin{equation}
n>\pi (n)\ln \pi (n)
\end{equation}
and
\begin{equation}
n<{\frac{\pi ^{-1}(n)}{\ln n}}.
\end{equation}
\end{lemma}

Using the function $\pi^{-1}(n)$ two basic entities of the present paper are
introduced: the sequence $\epsilon_{p_0}$ and the aggregate $r_{p_0}$.

\begin{definition}
The prime number sequence
\begin{equation}
\epsilon _{p_0}:p_0\in N,\ \ p_{k+1}=\pi ^{-1}(p_k),\ \ k=0,1,2,\ldots \ ,
\end{equation}
is called ''Eratosthenes progression of base $p_0$''. The infinite prime
number aggregate
\begin{equation}
r_{p_0}=\{p_{k+1}:p_0\in N,\ \ p_{k+1}=\pi ^{-1}(p_k),\ \ k=0,1,2,\ldots \}
\end{equation}
is called ''Eratosthenes ray of base $p_{0}$''.
\end{definition}

The following two auxiliary propositions are true for the Eratosthenes rays:

\begin{lemma}
Let two rays $r_{p_0^{\prime }}$ and $r_{p_0^{\prime \prime }}$ of bases
\begin{equation}
p_0^{\prime }<p_0^{\prime \prime }
\end{equation}
are given. Then if:

\begin{enumerate}
\item  $p_0^{\prime \prime }\in r_{p_0}^{\prime }$ the ray $r_{p_0^{\prime
\prime }}$ is contained in the ray $r_{p_0^{\prime }}$ ;

\item  $p_0^{\prime \prime }\overline{\in }r_{p_0^{\prime }}$ the ray $
r_{p_0^{\prime \prime }}$ does not contain common elements with the ray $
r_{p_0^{\prime }}$ and the inequality (12) implies the inequalities
\begin{equation}
p_k^{\prime }<p_k^{\prime \prime },\ \ k=0,1,2,\ldots \ .
\end{equation}
\end{enumerate}

\textbf{Proof}
\end{lemma}

\begin{enumerate}
\item  The inequality (12) gives the possibility to make the assumption that
$p_0^{\prime \prime }\in r_{p_0^{\prime }}$. That means that a number $
k^{*}>0$ exists so that $p_0^{\prime \prime }\equiv p_{k^{*}}^{\prime }$.
Now the inclusion $r_{p_0^{\prime \prime }}\subset r_{p_0^{\prime }}$
follows from the recursion $p_k^{\prime \prime }=\pi
^{-1}(p_{k+k^{*}}^{\prime }),\ \ k=0,1,2,\ldots $, given in (11).

\item  Let $p_0^{\prime \prime }\ \overline{\in }\ r_{p_0^{\prime }}$. Now
suppose a number $k^{*}>0$ exists so that $p_{k^{*}}^{\prime \prime }\in
r_{p_0^{\prime }}$. This implies the existence of a finite reverse sequence
\[
\ p_k^{\prime \prime }=\pi ^{-1}(p_{k+1}^{\prime \prime })\in r_{p_0^{\prime
}},\ \ \ where\ \ k=k^{*}-1,\ k^{*}-2,\ k^{*}-3,\ldots ,0,
\]
that terminates on the inclusion $p_0^{\prime \prime }\in r_{p_0}^{\prime }$
which contradicts the initial assumption $p_0^{\prime \prime }\ \overline{
\in }\ r_{p_0^{\prime }}$. Thus the ray $r_{p_0^{\prime \prime }}$ does not
have common elements with the ray $r_{p_0^{\prime }}$.

The inequality (12) gives the possibility to use the implication (7)
together with the recursion (11) that defines the elements of $
r_{p_0^{\prime }}$ and $r_{p_0^{\prime \prime }}$ and this leads to the
following chain of implications
\begin{eqnarray*}
p_0^{\prime }<p_0^{\prime \prime }&\Longrightarrow& p_1^{\prime }=\pi
^{-1}(p_0^{\prime })<p_1^{\prime \prime }=\pi ^{-1}(p_0^{\prime \prime }), \\
p_1^{\prime }<p_1^{\prime \prime }&\Longrightarrow& p_2^{\prime }=\pi
^{-1}(p_1^{\prime })<p_2^{\prime \prime }=\pi ^{-1}(p_1^{\prime \prime }), \\
 &\vdots&  \\
p_{k-1}^{\prime }<p_{k-1}^{\prime \prime }&\Longrightarrow& p_k^{\prime }=\pi
^{-1}(p_{k-1}^{\prime })<p_k^{\prime \prime }=\pi ^{-1}(p_k^{\prime \prime
}). \ \ [] \\
\end{eqnarray*}
\end{enumerate}

\begin{lemma}
The principal Eratosthenes rays do not intersect
\begin{equation}
\bigcap_{p_0\in \overline{C}}r_{p_0}=\emptyset
\end{equation}
where
\[
\ \overline{C}\equiv N\backslash P=\{1,4,6,8,9,10,12,\ldots \}.
\]
The union of all principal rays contains the union of all rays
\begin{equation}
\bigcup_{p_0\in P}r_{p_0}\subset \bigcup_{p_0\in \overline{C}}r_{p_0}.
\end{equation}
\textbf{Proof}
\end{lemma}

The proving method is by induction.

The first induction step starts considering the following 6 cases:

\begin{itemize}
\item  \emph{Case $p_0=1\in \overline{C}$}. The progression $\epsilon _1$
from (10) generates the first principal prime numbers ray $r_1$;

\item  \emph{Case $p_0=2\in r_1$ and $p_0=3\in r_1$.} The progressions $
\epsilon _2$ and $\epsilon _3$ from (10) determine the rays $r_2$ and $r_3$,
which according to lemma 2(1) do not contain any new elements compared to
the principal ray that is generated already. So the following inclusion is
true:
\begin{equation}
r_2\bigcup r_3\subset r_1;
\end{equation}

\item  \emph{Case $p_0=4\in \overline{C}$.} The progression $\epsilon _4$
from (10) generates the second principal prime numbers ray $r_4$, which
according to lemma 2(2) does not have common elements with the ray $r_1$.
Thus their intersection is void:
\begin{equation}
r_1\bigcap r_4=\emptyset ;
\end{equation}

\item  \emph{Case $p_0=5\in r_1$.} According to lemma 2(1) the Eratosthenes
progression $\epsilon _5$ generates the ray $r_5\subset r_1$. Together with
(16) this inclusion implies the new inclusion:
\begin{equation}
r_2\bigcup r_3\bigcup r_5\subset r_1;
\end{equation}

\item  \emph{Case $p_0=6\in \overline{C}$.} According to lemma 2(2) the
Eratosthenes progression $\epsilon _6$ generates the third principal prime
numbers ray, which does not contain common elements with the already
existing principal rays $r_1$ and $r_4$. Together with the relation (17)
this leads to a new particular case of the equality (14)
\[
\ r_1\bigcap r_4\bigcap r_6=\emptyset ;
\]

\item  \emph{Case $p_0=7\in r_4$.} According to lemma 2(1) the Eratosthenes
progression $\epsilon _7$ generates the ray $r_7\subset r_4$. This inclusion
together with the inclusion (18) imply the extended inclusion
\[
\ (r_2\bigcup r_3\bigcup r_5\bigcup r_7)\subset (r_1\bigcup r_4),
\]
On its turn this inclusion appears to be a particular case of the inclusion
(15).
\end{itemize}

The second induction step.
Suppose the relations
\begin{equation}
\bigcap_{p_0 \in C_n} r_{p_0}= \emptyset \ ,
\end{equation}
\begin{equation}
R_n \equiv \bigcup_{p_0 \in P_{c_n}} r_{p_0} \subset \bigcup_{p_0 \in C_n}
r_{p_0} \ ,
\end{equation}
where
\[
C_n = \bigcup_{i=1,2,3, \ldots , n} c_i,\ \ P_{c_n}=\{p \in P: \ p < c_n\},
\]
are satisfied up to the n-th arbitrary element $c_n > 7$ of the set
$\overline{C}$.

We shall prove that the relations (19) and (20) remain true for
$c_n = c_n+1$.
Let us consider the following two cases:

\begin{itemize}
\item  \emph{Case $p_0=c_n+1\in P$ .} We shall proceed starting from the
inequality
\begin{equation}
c_{n+1}>\pi (c_n+1)
\end{equation}
(the prime number $c_{n+1}$ is greater than its ordinal number $\pi
(c_n+1)$ ), which follows directly from inequality (8).

The right hand side of inequality (21) $\pi (c_n+1)$ is either a composite
or a prime number. Assume that $\pi (c_n+1)$ is a prime number. Then
applying the operation $\pi $ on both sides of (21) we obtain a double chain
of inequalities
\begin{equation}
c_{n+1}>\pi (c_n+1)>\pi ^2(c_n+1).
\end{equation}
Now we assume that the ''second number'' $\pi ^2(c_n+1)$ is again a prime
number and apply once more the operation $\pi $ on every element of the
double chain of inequalities (22). Thus we obtain a triple chain of
inequalities analogous to (22). The process can go further $\alpha $-times
obtaining $\alpha $-multiple chain of inequalities
\begin{equation}
c_{n+1}>\pi (c_n+1)>\pi ^2(c_n+1)>\ldots >\pi ^\alpha (c_n+1)
\end{equation}
until the ''$\alpha $-number'' $\pi ^\alpha (c_n+1)\equiv p_0^{*}\in
\overline{C}$.

Applying the operation $\pi $ on both sides of the recurrence formula (11)
for the ray $\epsilon _{p_0^{*}}$ and taking into account the identity (6)
we obtain the generator of the finite reverse sequence
\begin{equation}
p_i^{*}=\pi (p_{i+1}^{*}),\ \ i=\alpha ,\ \alpha -1,\ \alpha -2,\ \ldots ,0.
\end{equation}

The terms of the inequalities (23) are generated by the (24): $
c_n+1=p_\alpha ^{*}$, \ $\pi (c_n+1)=p_{\alpha -1}^{*},$ \ $\pi
^2(c^n+1)=p_{\alpha -2}^{*},$ \ \ldots , \ $\pi ^\alpha (c_n+1)=p_0^{*}$.
Thus the inequality $c_n+1>\pi ^\alpha (c_n+1)$ and the fact that $c_n+1$
belongs to the principal ray of base $\pi ^\alpha (c_n+1)$ (the latter is
shown by the sequence (24)) allows for the use of the lemma 2(1). This
implies the inclusion $r_{c_n+1}\subset R_n$, which together with eq. (20)
allows for the needed extension of the inclusion (20)
\[
\ r_{c_n+1}\bigcup \left( \bigcup_{p_0\in P_{c_n}}r_{p_0}\right)
=\bigcup_{p_0\in P_{c_n+1}}r_{p_0}\subset \bigcup_{p_0\in C_n}r_{p_0}\ .
\]

\item  \emph{Case $c_n+1=c_{n+1}\in \overline{C}$.} According to lemma 2(2)
the Eratosthenes progression $\epsilon _{c_{n+1}}$ generates the new
principal ray $r_{c_{n+1}}$, which is not contained in the union $R_n$. So
that due to the extension of relation (19)
\[
\ \bigcap_{p_0\in C_{n+1}}r_{p_0}=\emptyset \ .
\]
\end{itemize}

The third induction step takes into account the limits $C_n \longrightarrow
\overline{C}$ for $n \longrightarrow \infty$ and $P_{c_n +1} \longrightarrow
P$ for $n \longrightarrow \infty$. Thus we find out that the relations (19)
and (20) in that limit go to the relations (14) and (15) respectively. \ \
[] \\

\vspace{3mm}

Lemma 3 makes possible to show that the following basic theorem concerning
the separating of prime numbers into subsets with explicit law for the
determination of their elements actually takes place (Prime Number
Separating Theorem - PNST [6]).

\begin{theorem}
The set of prime numbers has a two-dimensional representation labeled by the
index $k$ and the base $p_0$
\begin{equation}
P=\bigcup_{p_0\in N\backslash P}\left\{ r_{p_0}=\{p_{k+1}:\ p_{k+1}=\pi
^{-1}(p_k),\ k=0,1,2,\ldots \}\right\} .
\end{equation}
\textbf{Proof}
\end{theorem}

The union of all Eratosthenes rays of bases covering the set of all natural
numbers
\begin{equation}
Q=\bigcup_{p_0\in N}\{p_{k+1}:p_{k+1}=\pi ^{-1}(p_k),\ \ \ k=0,1,2,\ldots \}
\end{equation}
can be represented in the form
\[
\ Q=Q_1\bigcup Q_2\ ,
\]
where
\[
\ Q_1=\{p_1:p_1=\pi ^{-1}(p_0),\ \ p_0\in N\},
\]
\[
\ Q_2=\bigcup_{p_1\in P}\{p_{k+1}:p_{k+1}=\pi ^{-1}(p_k),\ \ k=1,2,3,\ldots
\}.
\]
$Q_1$ coincides with $P$ according to the Euclid theorem (1). On the other
hand $Q_2$ is contained in $P$. Thus the right hand side of (26) coincides
with $P$ (i.e. $P\equiv Q$).

The set $Q$ can be also represented in the form
\begin{equation}
P \equiv Q = Q_3 \bigcup Q_4 \ ,
\end{equation}
where
\[
Q_3 = \bigcup_{p_0 \in N \backslash P} r_{p_0} \ \ and \ \ \ Q_4 =
\bigcup_{p_0 \in P} r_{p_0} \ .
\]
Owing to relation (15) from lemma 3 it follows that $Q_4 $ appears to be a
fraction of $Q_3$ ($Q_4 \subset Q_3$). Now it follows from (27) that $P
\equiv Q_3$. \ \ [] \\

\section{Implications of the prime number separating theorem}

Since $N\backslash P = C\bigcup \{1\} \equiv \overline{C}$ where $C$ is the
set of composite numbers it appears that PNST is mapping the elements of
the
extended set of composite numbers $p_0 \in \overline C$ into the infinite
prime number rays $r_{p_0}$. So that we come to
\begin{corollary}
There is an unique reciprocal mapping
\[
\ \varphi _1(p_0):\overline{C}\longrightarrow \ \ ^2P,
\]
which maps the elements $p_0$ of the set $\overline{C}$ into the principal
infinite prime number rays $r_{p_0}$.
\end{corollary}

Let us note that
\[
^{2}P = \{r_{p_0}\}_{p_0 \in N \backslash P} = \{p_{\mu \nu} \}_{
\begin{array}{l}
\mu = 1,2,3, \ldots \\
\nu = 1,2,3, \ldots
\end{array}
}
\]
denotes the prime number representation matrix introduced
by eq. (25).

The corollary 1 on its turn is showing that an analogous matrix
representation exists for the natural numbers too:
\[
\ ^2N=\{\overline{C},\ ^2P\}.
\]
The upper left hand side of the infinite matrix $^2N$ is given in the
Appendix.

Theorem 1 has diverse implications. Here we shall mention just one more. Let
us denote the following classes of Eratosthenes rays:
\begin{eqnarray*}
K_1 & = & r_1, \\
K_2 & = & \{r_i\}_{i=2^j} \ , \ j = 1,2,3, \ldots , \\
K_3 & = & \{r_i\}_{i = 3^j} \bigcup r_{2.3}, \ j = 1,2,3, \ldots , \\
K_5 & = & \{r_i \}_{i = 5^j} \bigcup r_{2.5} \bigcup r_{3.5}, \ j = 1,2,3,
\ldots , \\
& \vdots & \\
K_p & = & \{ r_i \}_{i = p^j} \bigcup \left( \bigcup_{\alpha \leq p , \alpha
\in P} r_{\alpha .p} \right), \ j = 1,2,3, \ldots \ .
\end{eqnarray*}

\begin{corollary}
There is an unique reciprocal mapping
\[
\ \varphi _2(p):\overline{P}\longrightarrow \{K_p\}_{p\in \overline{P}},
\]
which maps the elements of the set $\overline{P}=\{1\}\bigcup P$ into the
elements of the set of classes $\{K_p\}_{p\in \overline{P}}$. For the
equivalence classes $K_p$ the following set-theoretic equalities
\[
\ \bigcap_{p\in \overline{P}}K_p=\emptyset ,\ \ \ \bigcup_{p\in \overline{P}
}K_p=\ P.
\]
take place.
\end{corollary}

The following two propositions for the elements of the matrix $^{2}P$ rows
take place.

\begin{theorem}
The series
\[
\ s(p_0)=\sum_{i=1}^\infty {\frac 1{{p_i(p_0)}}},\ \ \ p_0\in N\backslash P
\]
are convergent. \newline
\textbf{Proof}
\end{theorem}

Following the inequalities (13) from lemma 2(2) the series $s(1)$
majorates all the series $s(p_0)$, for $p_0\in C$.  On its turn for
$i\geq 5$ the series $s(1)$ is majorated by the series $\sum_{i=1}^\infty
1/i^2$.  Indeed the inequalities (7), (9) and $5^2<31=p_5(1),$ imply the
following chain of inequalities \[ \ i^2<2i^2<\pi ^{-1}(i^2)<\pi
^{-1}(p_{i-1}(1)),\ \ for\ all\ i=6,7,8,\ldots \ \ .  \] \ \ [] \\

The estimate (8) from lemma 1 implies

\begin{lemma}
The spacing between two adjacent elements of a ray from $\{r_{p_0}\}_{p_0\in
\overline{C}}$ except for $p_{1,1}=2\in r_1$ only is estimated by
\begin{equation}
p_{(k+1)p_0}-p_{kp_0}>p_{kp_0}(\ln p_{kp_0}-1),\ \ \ for\ k=1,2,3,\ldots \ .
\end{equation}
\end{lemma}

Theorem 1 and its implications are applicable to all mathematical
constructions which include countable sets, or even sets containing finite
segments of $N$ or $P$ only.

\section{An analytic form of the prime number inner
distribution law}

The right hand side of relation (25) gives the definition for the recurrent
element of any ray $r_{p_0}, \ \ p_0 \in \overline{C}$ in terms of the
preceding element of the same ray only:
\begin{equation}
p_0 \in \overline{C}, \ \ p_{j+1} = \pi^{-1}(p_j), \ \ j=0,1,2,3, \ldots
 \ \ .
\end{equation}
This rule we comprehend to be the prime number inner distribution law ---
PNIDL.

The disadvantage of the law (29) is that the function $\pi^{-1}(x)$ (as well
as the function $\pi(x)$) can be realized by means of the Eratosthenes sieve
only. The Legendre formula (2) and its generalizations (see [8], p.343]) can
not be used for the purpose.

So far the attempts to derive analytic formulae for $\pi ^{-1}(p_j)$ for all
$j\times p_0$ led to formulae which do not allow for a new information on
prime numbers. All these formulae represent $\pi ^{-1}(x)$ as a discrete
function. Using these formulae it is only possible to determine these prime
numbers which are initially presupposed by their construction. These
formulae can not be extrapolated so as to increase the amount of prime
numbers (formula (2) is an example of that type of formulae).

We shall show here that the Eratosthenes rays can be approximated by
continuous functions which have extrapolation properties.

For an arbitrary row of the matrix $^2P$ from the Appendix we consider the
solution of the quadratic system of $m=2n$ equations
\begin{equation}
\left\{ \sum_{k=1}^n\alpha _{kp_0}q_{kp_0}(j)e^{-\beta _{kp_0}j}={\frac{{\ln
\ln \ln p_{1p_0}}}{{\ln \ln \ln p_{jp_0}}}}\right\} _{j=1,2,3,\ldots ,m}\ ,
\end{equation}
respective to the $n$ unknown pairs $\left\{ \alpha _{kp_0},\beta
_{kp_0}\right\} _{k=1,2,3,\ldots ,n}$, where
\begin{eqnarray*}
n &=&4\ \ in\ case\ of\ \ 1\leq p_0\leq 18\ \ (the\ first\ 11\ rays\ from\
Appendix) \\
n &=&3\ \ in\ case\ 20\leq p_0\leq 64\ \ (the\ succeeding\ 32\ rays\ from\
Appendix) \\
&&and \\
n &=&2\ \ in\ case\ of\ 65\leq p_0\leq 132\ \ (the\ last\ 57\ rays\ from\
Appendix)
\end{eqnarray*}
The polynomials $q_{kp_0}$ in system (30) take the values
\begin{eqnarray}
q_{kp_0}(j) &\equiv& 1 \ \ in\ case\ 1\leq p_0\leq 16,\ 20\leq p_0\leq 28,\
32\leq p_0\leq 132; \\
q_{1,18}(j) &\equiv &q_{2,18}(j)\equiv 1,\ \ q_{3,18}(j)\equiv
q_{4,18}(j)\equiv j\ \ for\ p_0=18; \\
&&and  \nonumber \\
q_{1,30}(j) &\equiv &q_{2,30}(j)\equiv 1,\ \ q_{3,30}(j)\equiv j,\ \ for\
p_0=30.
\end{eqnarray}

\begin{remark}
There are only two exceptions appearing for the first 100 rays from the
Appendix for which $q_{kp_0}(j)\neq 1$. These are $p_0=18$ and $p_0=30$.
\end{remark}

\begin{remark}
Constituting the system (30) for $p_0=1,4,6$ the first prime number $p_{1p_0}
$ is taken to be that for which the inequality $\ln \ln \ln (p_{1p_0})>0$
takes place for the first time. So that $p_{1,1}=31,\ \ p_{1,4}=17,\ \
p_{1,6}=41$.
\end{remark}

The principle feature of the system (30) is that it is exactly soluble. The
coefficients
\[
\left\{\alpha_{kp_0},\beta_{kp_0}\right\}_{
\begin{array}{l}
k = 1,2,3, \ldots , n \\
p_0 = 1,4,6, \ldots , 132
\end{array}
},
\]
satisfy the system (30) with residuals $<10^{-16}$. The solutions of the
system (30) have also the following properties:

\begin{enumerate}
\item  the amplitudes $\alpha _{kp}$ and the decrements $\beta _{kp}$ are
positive numbers;

\item  $\alpha _{kp}$ and $\beta _{kp}$ decrease when the index $k$ is
increasing;

\item  the increase of the amplitudes $\alpha _{kp_0}$ corresponds to an
increase of the decrements $\beta _{kp_0}$ (this trend is strictly
manifested for the first 26 rays);

\item  each pair of two consequent elements of the ray $r_{p_0}$ determine a
new term of the sum (30).
\end{enumerate}

The nonlinear systems (30) were analyzed by means of the program AFXY [9],
which determines the number of solutions and their accuracy according to the
method developed in [10] -- [14]. It was established that all systems (30)
-- (31) are uniquely soluble, while the systems (30) -- (32), (30) -- (33)
have triple solutions. In average the coefficients $\left\{ \alpha
_{kp_0},\beta _{kp_0}\right\} $ have 8 -- 10 correct decimal signs.

The solution of the system (30) led to the following approximation for the
function $\pi^{-1}(x)$ which covers all the elements contained in the first
100 Eratosthenes rays from the Appendix:
\begin{equation}
\tilde{\pi}^{-1}(x;p_0) = \exp{\exp{\exp{\frac{{\ln \ln \ln p_{1p_0}}}{{
\eta(x;p_0,n)}}}}}, \ \ \ x \in [1,m] \subset R^1,
\end{equation}
where
\[
\eta(x;p_0,n) = \sum_{k=1}^n \alpha_{kp_0}q_{kp_0}(x)e^{-\beta_{kp_0}x}.
\]

The formula
\begin{equation}
p_{jp_0} = \ \ round-off \ \ (\tilde{\pi}^{-1}(j;p_0))
\end{equation}
exactly reproduces the prime numbers from the Appendix up to the even number
$j^* \leq m$.

The function $\tilde{\pi}^{-1} (x)$ predicts the values of the new prime numbers $
p_{2n+1}^f$. As seen from Table where the mean accuracy values in respect to
$p_0$ are given
\[
\ \delta _n={\frac{{|p_{2n+1}-p_{2n+1}^f|_{100}}}{{p_{2n+1}}}}
\]
its prediction accuracy increases with the increase of the number n.

\textbf{Table} \newline

\begin{tabular}{|l||c|c|c|}
\hline
n & 2 & 3 & 4 \\ \hline
$\delta_n$ & 21\% -- 16\% & 5\% -- 1.7\% & 2.5\% -- 0.19\% \\ \hline
\end{tabular}

\vspace{7mm}

To compare with we point out that for $p_0 =4$ the relative accuracy of the
solution $p^f_{9,4}$ obtained from the equation

\begin{equation}
li \ x = p_{jp_0}
\end{equation}
when the right hand side is $p_{jp_0} = p_{8,4}$ equals to $0.004\%$

Let $l$ is the number of the correct decimal signs (the length of a computer
word) for which the arithmetic floating point operations in a given
computing environment (the computer) are produced. The above found
approximation for $\tilde \pi ^{-1}(x)$, which is limited with respect to
its domain of definition, suggests to check the hypothesis:

\begin{conjecture}
For any ray $r_{p_0},\ \ p_0\in \overline{C}$ there exist finite numbers $l$
and $n^{*}(l)$ such that the natural numbers $\tilde p_{mp_0}$ with $
m=2n^{*}(l)+1$ predicted by formulae (34), (35) lye closer to the prime
numbers $p_{mp_0}$ than the solutions $x$ of equation (36) when its right
hand sides are $p_{(m-1)p_0}$.
\end{conjecture}

\section{The set of origins of the Eratosthenes rays and
coagulates of prime numbers}

The set
\[
P_1 = \left\{ p_1(p_0) : p_1(p_0) = \pi^{-1}(p_0), \ \ p_0 \in \overline{C}
\right\} = \left\{p_{\mu 1} \right\}_{\mu = 1,2,3, \ldots}
\]
consists from the start-points (origins) of the principal rays
$\left\{r_{p_0}\right\}_{p_0 \in \overline{C}}$\ \ .

The set of prime numbers $^2P=\left\{ p_{\mu \nu }\right\} _{
\begin{array}{l}
\mu =1,2,3,\ldots  \\
\nu =1,2,3,\ldots
\end{array}
}$ is separated according to PNST onto two new infinite subsets
\[
\ ^2P=\left\{ P_1,P_{erat}\right\} ,\ \ where\ \ P_{erat}=\left\{ p_{\mu \nu
}\right\} _{
\begin{array}{l}
\mu =1,2,3,\ldots  \\
\nu =2,3,4,\ldots
\end{array}
} .
\]

The set $P_1$ takes a special place amidst the prime numbers. It seems
unlikely that a nonasymptotic distribution law similar to the PNIDL (29),
accounting for the rows of the matrix $P_{erat}$ too, can be derived for it.
Studying the nonasymptotic distribution of the elements of the subset $P_1$
amidst the prime numbers one can encounter its main peculiarity: $P_1$
contains all possible coagulates of prime numbers such as twin pairs and
other closely disposed sequences of prime numbers.

Let the sets of twin-pairs, twin-triples, twin-quadruples, twin-quintuples
of prime numbers greater than 5 are denoted respectively by
\begin{eqnarray*}
T_1 & = &\left\{ \overline{p}_i, \ \overline{p}_i +2 \right\}_{i=1,2,3,
\ldots , l_1}, \\
T_2 & = &\left\{\overline{p}_i, \ \overline{p}_i +2 , \ \overline{p}_i + 6
\right\}_{i=1,2,3, \ldots , l_2}, \\
T_3 & = & \left\{\overline{p}_i, \ \overline{p}_i +2, \ \overline{p}_i +6, \
\overline{p}_i+8 \right\}_{i=1,2,3, \ldots , l_3} , \\
T_4 & = & \left\{ \overline{p}_i, \ \overline{p}_i+2, \ \overline{p}_i +6, \
\overline{p}_i + 8, \ \overline{p}_i + 12 \right\}_{i=1,2,3, \ldots , l_4}.
\end{eqnarray*}
Following from lemma 4 is the basic proposition which binds the sets $T_1, \
T_2, \ T_3,$ and $T_4$ with the set $P_1$:

\begin{theorem}
None of the elements of the subsets $T_1,\ T_2,\ T_3$ and $T_4$ is not
contained in the Eratosthenes rays $\left\{ r_{p_0}\right\} _{p_0}\in
\overline{C}$. For any $i=1,2,3,\ldots ,l_k,\ (k=1,2,3,4)$ only one element
from pairs $T_1$, only two elements from the triplets $T_2$ and the
quadruplets $T_3$, and only three elements from the quintuplets $T_4$
could not be origins of certain principal ray
$\left\{ r_{p_0}\right\} _{p_0}\in \overline{C}$.
\end{theorem}

In theorem 3 the words ''could not be origins'' mean a probability whose
amount among the first $104683$ elements of the set $T_1$ does
not exceed the value $0.15$ ($1270$ twin-pairs are contained among these
first elements). The elements of the sets $T_2,\ T_3,$ and $T_4$ occur much
more rarely amidst the natural numbers compared to these of the set $T_1$.
Thus theorem 3 is showing that the essential part of the prime numbers
included in the sets $T_1,\ T_2,\ T_3$ and $T_4$ appear to be origins of
principal Eratosthenes rays.

A twin-pair assigns a characteristic ''arhythmicity of condensation'' among
the coagulates from $T_2,\ T_3$ and $T_4$ as well as among all elements of
the set $P_1$. Let us note that twin-pairs occur twice in $T_3$ and $T_4$.
For instance in the elements $t_{4,1}=\left\{ 11,13,17,19,23\right\} \in T_4$
and $t_{4,2}=\left\{ 101,103,107,109,113\right\} \in T_4$ these are the
pairs $\{11,13\},\{17,19\}$, and $\{101,103\},\ \{107,109\}$ respectively.

Theorem 1 and theorem 3 together show that the Eratosthenes rays "coagulate
between themselves" through the twin-pairs only and that according to the
character of the "ray coagulates" these twin-pairs are classified into two
new types:

\begin{enumerate}
\item  pairs such as $\{\overline{p}_{iu},\overline{p}_{iu}+2\}$ which
simultaneously give origin of two new rays;

\item  and the pairs $\{\overline{p}_{ib},\overline{p}_{ib}+2\}$ for which
one of the elements is the origin of a new ray while the other
''coagulates a new ray ''branching with an already existing ray.
\end{enumerate}

Thus
\begin{equation}
T_1 = T_{1u} \bigcup T_{1b} \ ,
\end{equation}
where
\[
T_{1u} = \{\overline{p}_{iu}, \overline{p}_{iu} +2 \}_{i = 1,2,3, \ldots ,
l_5} \subset P_1
\]
and
\[
T_{1b} = \{\overline{p}_{ib},\overline{p}_{ib} +2 \}_{i = 1,2,3, \ldots ,
l_6} = T_1 \backslash T_{2u}\ .
\]
As seen from the Appendix $\{71,73\}, \{101,103\}, \{137,139\}$ and $
\{149,151\}$ are examples of $u$-pairs while $\{11,13\}, \{17,19\}, \{59,61\}
$ and $\{107,109\}$ are examples of b-pairs.

The collective divisibility coefficient -- CDC
\begin{equation}
D(l,s)={\frac{{\left( d(l-s+1)+d(l-s+2)+\ldots +d(l)\right) }}{{s+2}}},
\end{equation}
where $l\geq 13,\ s\geq 3$ are natural numbers, $d(\lambda )$ is the number
of the divisors of the natural $\lambda $ except the unit, when applied to
the consequent triples ($s=3$) of natural numbers $p,\ p+1,\ p+2$ leads to a
new class of twin-pairs --- those for which $D(l,3)=1$ (the Strongly Related
Twins -- $T_{sr}$):
\begin{eqnarray} \nonumber
T_{sr} &=&\left\{
\{11,13\},\{17,19\},\{29,31\},\{41,43\},\{101,103\},\{137,139\},\ldots
\right\} = \\
&&\left\{ \overline{\overline{p}}_i,\overline{\overline{p}}_i+2\right\}
_{i=1,2,3,\ldots ,l_7}\ .
\end{eqnarray}
The even number $p+1$ from the twin-triple $\left\{ p,p+1,p+2\right\} $ is
factorized as follows
\begin{equation}
p+1=2.3.\sigma ,\ \ \ \sigma \in P.
\end{equation}
The numbers $\sigma =2,3,5,7,17,23$ from the factorization (40) correspond
to the twin-pairs (39) and form the origin of a new subset $\Sigma
\subset P$ which has a curious property: the last decimal digit of all $
\sigma \in \Sigma $ is either $3$ or $7$ (never $1$ or $9$; this property is
checked up to the $T_{sr}$--pair $\{47777,47779\}$).

\begin{remark}
\[
\ \min_{l\geq 13,s\geq 3}D(l,s)=D(l,3)=1\ .
\]
\end{remark}

\begin{remark}
The exceptions mentioned in Remark 1 are the even numbers 18 and 30
corresponed to $T_{sr}\ \ \{17,19\}$ and $\{29,31\}$.
\end{remark}

The decomposition (39) is carried over the elements of the set $T_{sr}$ too
\[
\ T_{sr}=T_{(sr)u}\bigcup T_{(sr)b}\ ,
\]
where
\[
\ T_{(sr)u}=\left\{ \overline{\overline{p}}_{iu},\ \overline{\overline{p}}
_{iu}+2\right\} _{i=1,2,3,\ldots ,l_8}
\]
and
\[
\ T_{(sr)b}=\left\{ \overline{\overline{p}}_{ib},\ \overline{\overline{p}}
_{ib}+2\right\} _{i=1,2,3,\ldots ,l_9}.
\]
The existence of a potential infinity $l_k=\infty $ for the lengths $l_k,\
k=1,2,3,\ldots ,9$ of the introduced sets of prime numbers is still not
proved even for the simplest case $k=1$.

The following extension of the V. Brun [15] theorem is possible:

\begin{theorem}
If $l_k=\infty $ \ for $k=1,5,6,7,8,9$ the series
\[
\ \sum_{i=1}^\infty \left( {\frac 1{{v_i}}}+{\frac 1{{v_i+2}}}\right) ,\ \ \
where\ v_i=\overline{p}_i,\ \overline{p}_{iu},\ \overline{p}_{ib},\
\overline{\overline{p}}_i,\ \overline{\overline{p}}_{iu},\ \overline{
\overline{p}}_{ib}
\]
are convergent.
\end{theorem}

In case of $l_2=l_3=l_4=\infty $ a statement analogous to theorem 4 will
take place for the elements of the sets $T_2,T_3,$ and $T_4$.

Besides the sets $T_2$ and $T_3$ one can consider the sets of coagulates of
the type $\{p_i, p_i+4, p_i+6\}, \ \{p_i, p_i+2, p_i+8\}, \ \{p_i, p_i+4,
p_i+6, p_i+10\}$ and $\{p_i, p_i+2, p_i+8, p_i+12\}$ too. For these sets the
analogues of theorems 3 and 4 are also correct.

Let us introduce the set
\[
\ \tilde P_1=P_1\backslash \left( T_1\bigcup T_2\bigcup T_3\bigcup
T_4\right) \equiv \left\{ \tilde p_{i1}\right\} _{i=1,2,3,\ldots }.
\]
Treating numerically the origins of the partial sums of the reciprocal
values of the elements of the columns of the matrix $^2P$ we come to a
hypothesis concerning the set of origins of the principal rays of $P_1$
which is opposite in since of the V. Brun theorem:

\begin{conjecture}
The series
\begin{equation}
\sum_{i=1}^\infty {\frac 1{{w_i}}},\ \ \ where\ w_i=p_{i1},\ \tilde p_{i1},\
\ p_{i\nu }\ \ \ (\nu =2,3,4,\ldots )
\end{equation}
is divergent.
\end{conjecture}

\begin{remark}
The divergency of the series (41) is very slow for $w_i=p_{i1},\ \tilde
p_{i1}$ and even more slow for $w_i=p_{i\nu }$. What is more the slowness is
growing with $\nu \rightarrow \infty $ (?).
\end{remark}

The coagulates of primes elements of the sets $T_1, T_2, T_3$ and $T_4$
considered so far should be generalized as
\[
coag \left(p,\left\{u_i\right\}_{i=1,2,3, \ldots, l } \right) =
\left\{p, p +u_1, p+u_2, \ldots , p+u_l \right\}
\]
where the steps $\{u_i\}$ should be
smaller than the steps $p_{(k+1)p_0} - p_{kp_0} = \Delta\epsilon_{kp_0}$ in
the nearest rays while the rays and their steps themselves (i.e. the
particular $p_0$ and $k$) are determined by the initial prime number $p$ and
the length $l$ of the generalized coagulate. The pointed out relativity of
the values of $p, l$ and $\Delta\epsilon_{kp_0}$ in $coag\left\{p,\left\{u_i
\right\}_{i=1,2, \ldots ,l} \right\}$ is used in order to obtain the
conditions under which the set of generalized \textsl{coags} satisfy
theorems analogous to theorems 3 and 4.
\\\\
The general meaning of theorems 3 and 4, and
conjecture 2 is that the main information on the originality of the
nonasymptotic prime-numbers distribution is concentrated in the set of
start--points $P_1$; by its nature the set $P_1$ contains ''enough
amount'' of prime-numbers and they are disposed ''enough closely''
(conjecture 2), however only the set of the type of $T_1,T_2,T_3$ and
$T_4$, which contain ''relatively small amount of elements'' (theorem 4)
introduce the characteristic, unique arrhythmia of concentration of
elements of the set $P_1$.  \newpage

\
section*{Appendix}

\textbf{The left upper corner of the infinite matrix ${}^2N$.} %\newline

\begin{tabular}{|r|rrrrr|}
\hline
1 & 2 & 3 & 5 & 11 & 31 \\

& 127 & 709 & 5381 & 52711 & 648391 \\
& 9737333 & 174440041... &  &  &  \\ \hline
4 & 7 & 17 & 59 & 277 & 1787 \\
& 15299 & 167449 & 2269733 & 37139213 & 718064159... \\ \hline
6 & 13 & 41 & 179 & 1063 & 8527 \\
& 87803 & 1128889 & 17624813 & 326851121... &  \\ \hline
8 & 19 & 67 & 331 & 2221 & 19577 \\
& 219613 & 3042161 & 50728129... &  &  \\ \hline
9 & 23 & 83 & 431 & 3001 & 27457 \\
& 318211 & 4535189 & 77557187... &  &  \\ \hline
10 & 29 & 109 & 599 & 4397 & 42043 \\
& 506683 & 7474967 & 131807699... &  &  \\ \hline
12 & 37 & 157 & 919 & 7193 & 72727 \\
& 919913 & 14161729 & 259336153... &  &  \\ \hline
14 & 43 & 191 & 1153 & 9319 & 96797 \\
& 1254739 & 19734581 & 368345293... &  &  \\ \hline
15 & 47 & 211 & 1297 & 10631 & 112129 \\
& 1471343 & 23391799 & 440817757... &  &  \\ \hline
16 & 53 & 241 & 1523 & 12763 & 137077 \\
& 1828669 & 29499439 & 563167303... &  &  \\ \hline
18 & 61 & 283 & 1847 & 15823 & 173867 \\
& 2364361 & 38790341 & 751783477... &  &  \\ \hline
20 & 71 & 353 & 2381 & 21179 & 239489 \\
& 3338989 & 56011909... &  &  &  \\ \hline
21 & 73 & 367 & 2477 & 22093 & 250751 \\
& 3509299 & 59053067... &  &  &  \\ \hline
22 & 79 & 401 & 2749 & 24859 & 285191 \\
& 4030889 & 68425619... &  &  &  \\ \hline
24 & 89 & 461 & 3259 & 30133 & 352007 \\
& 5054303 & 87019979... &  &  &  \\ \hline
25 & 97 & 509 & 3637 & 33967 & 401519 \\
& 5823667 & 101146501... &  &  &  \\ \hline
26 & 101 & 547 & 3943 & 37217 & 443419 \\
& 6478961 & 113256643... &  &  &  \\ \hline
\end{tabular}

\newpage

\textbf{The left upper \ldots Continuation} %\newline

\begin{tabular}{|r|rrrrr|}
\hline
27 & 103 & 563 & 4091 & 38833 & 464939 \\
& 6816631 & 119535373... &  &  &  \\ \hline
28 & 107 & 587 & 4273 & 40819 & 490643 \\
& 7220981... &  &  &  &  \\ \hline
30 & 113 & 617 & 4549 & 43651 & 527623 \\
& 7807321... &  &  &  &  \\ \hline
32 & 131 & 739 & 5623 & 55351 & 683873 \\
& 10311439... &  &  &  &  \\ \hline
33 & 137 & 773 & 5869 & 57943 & 718807 \\
& 10875147... &  &  &  &  \\ \hline
34 & 139 & 797 & 6113 & 60647 & 755387 \\
& 11469013... &  &  &  &  \\ \hline
35 & 149 & 859 & 6661 & 66851 & 839483 \\
& 12838937... &  &  &  &  \\ \hline
36 & 151 & 877 & 6823 & 68639 & 864013 \\
& 13243033... &  &  &  &  \\ \hline
38 & 163 & 967 & 7607 & 77431 & 985151 \\
& 15239333... &  &  &  &  \\ \hline
39 & 167 & 991 & 7841 & 80071 & 1021271 \\
& 15837299... &  &  &  &  \\ \hline
40 & 173 & 1031 & 8221 & 84347 & 1080923 \\
& 16827557... &  &  &  &  \\ \hline
42 & 181 & 1087 & 8719 & 90023 & 1159901 \\
& 18143603... &  &  &  &  \\ \hline
44 & 193 & 1171 & 9461 & 98519 & 1278779 \\
& 20137253... &  &  &  &  \\ \hline
45 & 197 & 1201 & 9739 & 101701 & 1323503 \\
& 20890789... &  &  &  &  \\ \hline
46 & 199 & 1217 & 9859 & 103069 & 1342907 \\
& 21219089... &  &  &  &  \\ \hline
48 & 223 & 1409 & 11743 & 125113 & 1656649 \\
& 26548261... &  &  &  &  \\ \hline
49 & 227 & 1433 & 11953 & 127643 & 1693031 \\
& 27170047... &  &  &  &  \\ \hline
50 & 229 & 1447 & 12097 & 129229 & 1715761 \\
& 27560453... &  &  &  &  \\ \hline
51 & 233 & 1471 & 12301 & 131707 & 1751411 \\
& 28171007... &  &  &  &  \\ \hline
\end{tabular}

\newpage

\textbf{The left upper \ldots Continuation} %\newline

\begin{tabular}{|r|rrrrr|}
\hline
52 & 239 & 1499 & 12547 & 134597 & 1793237 \\
& 28889363... &  &  &  &  \\ \hline
54 & 251 & 1597 & 13469 & 145547 & 1950629 \\
& 31599859... &  &  &  &  \\ \hline
55 & 257 & 1621 & 13709 & 148439 & 1993039 \\
& 32332763... &  &  &  &  \\ \hline
56 & 263 & 1669 & 14177 & 153877 & 2071583 \\
& 33691309... &  &  &  &  \\ \hline
57 & 269 & 1723 & 14723 & 160483 & 2167937 \\
& 35368547... &  &  &  &  \\ \hline
58 & 271 & 1741 & 14867 & 162257 & 2193689 \\
& 35815873... &  &  &  &  \\ \hline
60 & 281 & 1823 & 15641 & 171697 & 2332537 \\
& 38235377... &  &  &  &  \\ \hline
62 & 293 & 1913 & 16519 & 182261 & 2487943... \\ \hline
63 & 307 & 2027 & 17627 & 195677 & 2685911... \\ \hline
64 & 311 & 2063 & 17987 & 200017 & 2750357... \\ \hline
65 & 313 & 2081 & 18149 & 202001 & 2779781... \\ \hline
66 & 317 & 2099 & 18311 & 204067 & 2810191... \\ \hline
68 & 337 & 2269 & 20063 & 225503 & 3129913... \\ \hline
69 & 347 & 2341 & 20773 & 234293 & 3260657... \\ \hline
70 & 349 & 2351 & 20899 & 235891 & 3284657... \\ \hline
72 & 359 & 2417 & 21529 & 243781 & 3403457... \\ \hline
74 & 373 & 2549 & 22811 & 259657 & 3643579... \\ \hline
75 & 379 & 2609 & 23431 & 267439 & 3760921... \\ \hline
76 & 383 & 2647 & 23801 & 271939 & 3829223... \\ \hline
77 & 389 & 2683 & 24107 & 275837 & 3888551... \\ \hline
78 & 397 & 2719 & 24509 & 280913 & 3965483... \\ \hline
80 & 409 & 2803 & 25423 & 292489 & 4142053... \\ \hline
81 & 419 & 2897 & 26371 & 304553 & 4326473... \\ \hline
82 & 421 & 2909 & 26489 & 305999 & 4348681... \\ \hline
84 & 433 & 3019 & 27689 & 321017 & 4578163... \\ \hline
85 & 439 & 3067 & 28109 & 326203 & 4658099... \\ \hline
86 & 443 & 3109 & 28573 & 332099 & 4748047... \\ \hline
87 & 449 & 3169 & 29153 & 339601 & 4863959... \\ \hline
88 & 457 & 3229 & 29803 & 347849 & 4989697... \\ \hline
90 & 463 & 3299 & 30557 & 357473 & 5138719... \\ \hline
\end{tabular}

\newpage

\textbf{The left upper \ldots End.} %\newline

\begin{tabular}{|r|rrrrr|}
\hline
91 & 467 & 3319 & 30781 & 360293 & 5182717... \\ \hline
92 & 479 & 3407 & 31667 & 371981 & 5363167... \\ \hline
93 & 487 & 3469 & 32341 & 380557 & 5496349... \\ \hline
94 & 491 & 3517 & 32797 & 386401 & 5587537... \\ \hline
95 & 499 & 3559 & 33203 & 391711 & 5670851... \\ \hline
96 & 503 & 3593 & 33569 & 396269 & 5741453... \\ \hline
98 & 521 & 3733 & 35023 & 415253 & 6037513... \\ \hline
99 & 523 & 3761 & 35311 & 418961 & 6095731... \\ \hline
100 & 541 & 3911 & 36887 & 439357 & 6415081... \\ \hline
102 & 557 & 4027 & 38153 & 455849 & 6673993... \\ \hline
104 & 569 & 4133 & 39239 & 470207 & 6898807... \\ \hline
105 & 571 & 4153 & 39451 & 472837 & 6940103... \\ \hline
106 & 577 & 4217 & 40151 & 481847 & 7081709... \\ \hline
108 & 593 & 4339 & 41491 & 499403 & 7359427... \\ \hline
110 & 601 & 4421 & 42293 & 510031 & 7528669... \\ \hline
111 & 607 & 4463 & 42697 & 515401 & 7612799... \\ \hline
112 & 613 & 4517 & 43283 & 522829 & 7730539... \\ \hline
114 & 619 & 4567 & 43889 & 530773 & 7856939... \\ \hline
115 & 631 & 4663 & 44879 & 543967 & 8066533... \\ \hline
116 & 641 & 4759 & 45971 & 558643 & 8300687... \\ \hline
117 & 643 & 4787 & 46279 & 562711 & 8365481... \\ \hline
118 & 647 & 4801 & 46451 & 565069 & 8402833... \\ \hline
119 & 653 & 4877 & 47297 & 576203 & 8580151... \\ \hline
120 & 659 & 4933 & 47857 & 583523 & 8696917... \\ \hline
121 & 661 & 4943 & 47963 & 584999 & 8720227... \\ \hline
122 & 673 & 5021 & 48821 & 596243 & 8900383... \\ \hline
123 & 677 & 5059 & 49207 & 601397 & 8982923... \\ \hline
124 & 683 & 5107 & 49739 & 608459 & 9096533... \\ \hline
125 & 691 & 5189 & 50591 & 619739 & 9276991... \\ \hline
126 & 701 & 5281 & 51599 & 633467 & 9498161... \\ \hline
128 & 719 & 5441 & 53353 & 657121 & 9878657... \\ \hline
129 & 727 & 5503 & 54013 & 665843... &  \\ \hline
130 & 733 & 5557 & 54601 & 673793... &  \\ \hline
132 & 743 & 5651 & 55681 & 688249... &  \\ \hline
. & . & . & . & . \ldots &  \\ \hline
. & . & . & . & . \ldots &  \\ \hline
. & . & . & . & . \ldots &  \\ \hline
\end{tabular}

\newpage
\section*{References}
\hspace*{0.15mm} \indent 1) Ch.~J. de la Vallee-Poussin,\textit{Recherhe
analytiues sur la theorie de nombers premiers}, Ann, Soc Sci., Bruxelles,
$\underline{20}$, (1899).

2) H.G. Diamond, \textit{Elementary methods in the study of the distribution
of prime numbers}, Bull Amer. Math. Soc., $\underline{7}$, (1982), 553 --
589.

3) D.J. Newman, \textit{Simple analytic proof of the prime number theorem},
Amer. Math. Monthly, $\underline{87}$, (1980), 693 -- 696.

4) A.F. Lavrik, \textit{On the distribution of primes based on
I.M.Vinogradov's method of trigonometric sums}, Trudi Mat. Inst. Steklov,
$\underline{64}$, (1961), 90 -- 125 (in Russian).

5) E.C. Titchmarsh, \textit{The theory of Riemann zeta-function}, Clarendon
Press, (1951).

6) L. Alexandrov, \textit{Multiple Eratosthenes sieve and the prime number
distribution on the plane}, Sofia (1964) (unpublished).

7) B. Rosser, \textit{The n-th prime is greater than $n\log n$}, Proc.
London Math. Soc., $\underline{45}$, $N^0 2$, (1939), 2144.

8) A. Bukhstab, \textit{Number theory}, Moscow (1960) (in Russian).

9) L. Alexandrov, Alex Karailiev, \textit{Program AFXY (Analyze
$Fx=y$)}, Menlo Park, CA, (1998) (private program system).

10) L. Alexandrov, \textit{Regularized Newton--Kantorovich computational
processes}, J. of Comp. Math. and Math. Phys., $\underline{2}$, $N^0 1$,
(1971), 36 -- 41, (in Russian).

11) L. Alexandrov, V. Gadjokov, \textit{Analysis of latent regularities by
means of regularized iteration processes}, J. of Radioanal. Chemistry,
$\underline{9}$, (1971), 279 -- 292.

12) L. Alexandrov, \textit{Regularized trajectories of the Newton--type
approximation for solving non-linear equations},J of Differential
Equations, $\underline{13}$, $N^0 7$, (1977), 1281 -- 1292, (in Russian).

13) L. Alexandrov, M. Drenska, D. Karadjov, P. Morozov, \textit{Inverse
relativistic problem for the bound quark--antiquark states}, J. of
Computational Physics, $\underline{45}$, $N^0 2$, (1982), 291 -- 299.

14) L. Alexandrov, M. Drenska, D. Karadjov, \textit{Program REGN for
solving numerically non-linear systems via regularized Gauss--Newton
methods}, RSIC/PSR--165, ORNL, Oak Ridge, Tennessee, (1984).

15) V. Brun, \textit{La serie 1/5 + 1/7 + 1/11 + \ldots ou le denominateurs
sont nombre premiers jumeaux et convergente ou finie}`, Bull. Sci. Math., $
\underline{43}$, $N^0 2$, (1919) 100 -- 104; 124 -- 128.
\end{document}